\documentclass[11pt]{amsart}

\usepackage[T1]{fontenc}
\usepackage[utf8]{inputenc}
\usepackage{lmodern}
\usepackage{amsmath,amssymb,amsthm}
\usepackage[margin=1in]{geometry}
\usepackage[hidelinks]{hyperref}

\numberwithin{equation}{section}
\newtheorem{theorem}{Theorem}
\newtheorem{corollary}{Corollary}
\newtheorem{lemma}{Lemma}
\theoremstyle{remark}

\DeclareMathOperator*{\esssup}{ess\,sup}
\newcommand{\T}{\mathbb T}
\newcommand{\BMO}{\mathrm{BMO}}

\begin{document}
\title[A BMO Extension of Oskolkov's Criterion]{A BMO Extension of Oskolkov's Criterion for Rectangular
Trigonometric Fourier Sums}
\author{Ushangi Goginava}
\address{Department of Mathematical Sciences, United Arab Emirates
University, P.O. Box 15551, Al Ain, Abu Dhabi, United Arab Emirates}
\email{zazagoginava@gmail.com}
\email{ugoginava@uaeu.ac.ae}
\date{}
\subjclass[2020]{Primary 42B05; Secondary 42A20, 42B35}
\keywords{Rectangular Fourier sums, almost-everywhere convergence,
sectional BMO, modulus of continuity}

\begin{abstract}
We give a sufficient condition in terms of sectional BMO for
almost-everywhere Pringsheim convergence of rectangular trigonometric
Fourier sums. For continuous functions, the condition includes
Oskolkov's criterion, and an example shows that the inclusion is
strict.
\end{abstract}
\maketitle

\section{Introduction}

Let \(\T=\mathbb R/\mathbb Z\), with normalized Lebesgue measure \(m\),
and write \(m_2=m\otimes m\) for the product measure on \(\T^2\).
For \(f\in L^1(\T^2)\), put
\begin{align}
\widehat f(k,\ell)
  &=\int_{\T^2}f(x,y)e^{-2\pi i(kx+\ell y)}\,dx\,dy, \\
S_{M,N}f(x,y)
  &=\sum_{|k|\le M}\sum_{|\ell|\le N}
    \widehat f(k,\ell)e^{2\pi i(kx+\ell y)}.
    \label{eq:rectangular-sum}
\end{align}
The Pringsheim limit means that \(\min(M,N)\to\infty\), with the two
indices otherwise independent.  For \(f\in C(\T^2)\), define
\[
\omega(f,\delta)_C
 =\sup_{\substack{s^2+t^2\le\delta^2\\(x,y)\in\T^2}}
       |f(x+s,y+t)-f(x,y)|.
\]

Fefferman \cite{Fefferman1971} constructed a continuous function whose
rectangular Fourier sums diverge in the Pringsheim sense.  Bakhbukh and
Nikishin \cite{BN} obtained a divergence example with
\begin{equation}
\omega(f,\delta)_C
 =O\!\left(\frac{1}{\log(1/\delta)}\right)
 \qquad(\delta\downarrow0).
\label{eq:BN}
\end{equation}
Oskolkov \cite[Theorem 5]{Oskolkov1974} proved almost-everywhere
Pringsheim convergence under the condition
\begin{equation}
\omega(f,\delta)_C
 =o\!\left(\frac{1}
 {\log(1/\delta)\log\log\log(1/\delta)}\right).
\label{eq:Oskolkov}
\end{equation}
Masyutin \cite{Masyutin2026} later obtained criteria for \(L^p\)
functions, \(1<p<\infty\), in terms of their \(L^p\) modulus of
continuity.

For \(u\in L^1(\T)\), write
\begin{equation}
[u]_{\BMO(\T)}
 =\sup_{I\subset\T}\frac1{|I|}\int_I|u-u_I|\,dx,
 \qquad u_I=\frac1{|I|}\int_Iu\,dx,
\label{eq:bmo}
\end{equation}
where \(I\) ranges over arcs.  Define the sectional seminorm
\begin{equation}
[f]_{\mathrm{sBMO}}
=\max\!\left\{
 \esssup_{y\in\T}[f(\cdot,y)]_{\BMO(\T)},
 \esssup_{x\in\T}[f(x,\cdot)]_{\BMO(\T)}
 \right\}.
\label{eq:sectional-bmo}
\end{equation}
We write \(f\in\mathrm{sBMO}(\T^2)\) when \(f\in L^1(\T^2)\) and
\eqref{eq:sectional-bmo} is finite.
For \(0<\delta\le 1/2\), put
\begin{align}
\Omega_1(f,\delta)
 &=\sup_{|h|\le\delta}\esssup_{y\in\T}
   [f(\cdot+h,y)-f(\cdot,y)]_{\BMO(\T)},\\
\Omega_2(f,\delta)
 &=\sup_{|h|\le\delta}\esssup_{x\in\T}
   [f(x,\cdot+h)-f(x,\cdot)]_{\BMO(\T)},\\
\Omega(f,\delta)&=\max\{\Omega_1(f,\delta),\Omega_2(f,\delta)\}.
\label{eq:modulus}
\end{align}
All asymptotic conditions below concern \(\delta\downarrow0\).

\begin{theorem}\label{thm:main}
Let \(f\in\mathrm{sBMO}(\T^2)\), \(n_q=2^{2^q}\), and
\[
a_q=2^q\Omega\!\left(f,\frac1{n_q+1}\right).
\]
If, for every \(\eta>0\),
\begin{equation}
\sum_{q=1}^{\infty}\exp(-\eta/a_q)<\infty,
\label{eq:new-condition}
\end{equation}
where a summand with \(a_q=0\) is interpreted as zero, then
\(S_{M,N}f(x,y)\to f(x,y)\) for almost every \((x,y)\) as
\(\min(M,N)\to\infty\).
\end{theorem}

\begin{corollary}\label{cor:osk-bmo}
Let \(f\in\mathrm{sBMO}(\T^2)\).  If
\[
\Omega(f,\delta)
=o\!\left(\frac1{\log(1/\delta)
  \log\log\log(1/\delta)}\right),
\]
then \(S_{M,N}f\to f\) almost everywhere in the Pringsheim sense.
\end{corollary}

\begin{proof}
For \(\delta_q=(2^{2^q}+1)^{-1}\), we have
\(\log(1/\delta_q)\asymp 2^q\) and
\(\log\log\log(1/\delta_q)\sim\log q\).  Consequently
\(a_q=o(1/\log q)\).  Given \(\eta>0\), eventually
\(\eta/a_q\ge2\log q\) whenever \(a_q>0\); thus the summands in
\eqref{eq:new-condition} are at most \(q^{-2}\).  Apply
Theorem~\ref{thm:main}.
\end{proof}

For \(f\in C(\T^2)\), the bound
\begin{equation}
\Omega(f,\delta)\le2\omega(f,\delta)_C
\label{eq:continuous-domination}
\end{equation}
shows that Corollary~\ref{cor:osk-bmo} includes Oskolkov's criterion.
The next example shows that \eqref{eq:new-condition} can hold for a
continuous function that does not satisfy Oskolkov's hypothesis.

\begin{theorem}\label{thm:example}
There is \(F\in C(\T^2)\) such that, for every \(\eta>0\),
\begin{equation}
\sum_{q=1}^{\infty}
 \exp\!\left(-\frac{\eta}
 {2^q\omega(F,(n_q+1)^{-1})_C}\right)<\infty,
\label{eq:example-condition}
\end{equation}
with the usual zero-summand convention, but
\begin{equation}
r\log\log r\;\omega(F,2^{-r})_C\not\longrightarrow0
\qquad(r\to\infty).
\label{eq:example-failure}
\end{equation}
Thus Oskolkov's hypothesis need not hold when
Theorem~\ref{thm:main} applies.
\end{theorem}

\section{One-dimensional estimates}

Write \(S_n^{(1)}\) and \(S_n^{(2)}\) for Fourier partial sums in \(x\)
and \(y\), respectively.  Explicitly,
\begin{align*}
S_n^{(1)}f(x,y)
 &=\sum_{|k|\le n}\!
   \left(\int_\T f(s,y)e^{-2\pi iks}\,ds\right)e^{2\pi ikx},\\
S_n^{(2)}f(x,y)
 &=\sum_{|\ell|\le n}\!
   \left(\int_\T f(x,t)e^{-2\pi i\ell t}\,dt\right)e^{2\pi i\ell y}.
\end{align*}
By Fubini,
\begin{equation}
S_{M,N}=S_M^{(1)}S_N^{(2)}=S_N^{(2)}S_M^{(1)}.
\label{eq:factor}
\end{equation}
A polynomial in one coordinate may have measurable coefficients
depending on the other coordinate.

\begin{lemma}
\label{lem:jackson}
There is an absolute \(C_J>0\) such that, for
\(f\in\mathrm{sBMO}(\T^2)\), \(n\ge1\), and \(\nu=1,2\), there exists
a polynomial \(P_n^{(\nu)}\) of degree at most \(n\) in coordinate
\(\nu\).  The residual \(g_n^{(\nu)}=f-P_n^{(\nu)}\) has mean zero in
the active coordinate, and
\begin{align}
\esssup_y[g_n^{(1)}(\cdot,y)]_{\BMO(\T)}
 &\le C_J\Omega_1\!\left(f,\frac1{n+1}\right),
 \label{eq:jackson1}\\
\esssup_x[g_n^{(2)}(x,\cdot)]_{\BMO(\T)}
 &\le C_J\Omega_2\!\left(f,\frac1{n+1}\right).
 \label{eq:jackson2}
\end{align}
\end{lemma}

\begin{proof}
Set \(m=\lfloor n/2\rfloor+1\) and use
\[
F_m(t)=\frac1m\left(\frac{\sin(\pi m t)}{\sin(\pi t)}\right)^2,
\qquad J_n(t)=\frac{F_m(t)^2}{\int_\T F_m(s)^2\,ds}.
\]
Thus \(F_m\) has degree \(m-1\), and the nonnegative kernel
\(J_n\) has degree \(2(m-1)\le n\) and integral one.  The bounds
\(F_m(t)\le C\min\{m,(m|t|^2)^{-1}\}\) and
\(\int_\T F_m^2\asymp m\) give
\begin{equation}
\int_\T |t|J_n(t)\,dt\le\frac{C}{n+1},
\label{eq:first-moment}
\end{equation}
where \(t\in[-1/2,1/2]\) is the representative of a point of \(\T\).
For \(\nu=1\), take
\[
P_n^{(1)}(x,y)=\int_\T f(x-t,y)J_n(t)\,dt.
\]
This is a polynomial in \(x\) with measurable coefficients in \(y\).
  With
\(\delta=(n+1)^{-1}\), Minkowski's inequality for the BMO seminorm
gives
\[
[g_n^{(1)}(\cdot,y)]_{\BMO}
 \le\int_\T J_n(t)
  [f(\cdot,y)-f(\cdot-t,y)]_{\BMO}\,dt.
\]
For a fixed \(t\in[-1/2,1/2]\), set
\[
N(t)=\max\!\left\{1,\left\lceil\frac{|t|}{\delta}\right\rceil\right\},
\qquad h=\frac{t}{N(t)}.
\]
Then \(|h|\le\delta\) and \(N(t)\le1+|t|/\delta\).
Write \(u_y=f(\cdot,y)\) and \((\tau_su)(x)=u(x+s)\).
The identity
\[
u_y-\tau_{-t}u_y
 =\sum_{j=0}^{N(t)-1}
   \tau_{-jh}(u_y-\tau_{-h}u_y)
\]
is a telescoping sum.  For each fixed \(t\), the definition of
\(\Omega_1\), together with \(|h|\le\delta\), gives
\([u_y-\tau_{-h}u_y]_{\BMO}\le\Omega_1(f,\delta)\)
for almost every \(y\).  The triangle inequality and translation
invariance of the BMO seminorm therefore imply
\[
[f(\cdot,y)-f(\cdot-t,y)]_{\BMO}
 \le N(t)\Omega_1(f,\delta)
 \le\left(1+\frac{|t|}{\delta}\right)\Omega_1(f,\delta).
\]
The exceptional set of \(y\)'s may depend on \(t\); Fubini's
theorem shows that this estimate holds for almost
every \(t\) for almost every \(y\).  We may therefore integrate it
against \(J_n(t)\).  Since \(J_n\ge0\), \(\int_\T J_n=1\), and
\eqref{eq:first-moment} gives
\(\delta^{-1}\int_\T|t|J_n(t)\,dt\le C\), we obtain
\[
[g_n^{(1)}(\cdot,y)]_{\BMO}
 \le\Omega_1(f,\delta)
 \left(1+\frac1\delta\int_\T|t|J_n(t)\,dt\right)
 \le(1+C)\Omega_1(f,\delta)
\]
for almost every \(y\).  Taking the essential supremum proves
\eqref{eq:jackson1}; interchanging \(x\) and \(y\) proves
\eqref{eq:jackson2}.
\end{proof}

For \(u\in L^2(\T)\), let
\[
\mathcal Hu
  =\sum_{k\in\mathbb Z}(-i\,\operatorname{sgn}k)
          \widehat u(k)e^{2\pi ikx},
\qquad
\widetilde S_nu=S_n(\mathcal Hu).
\]
We use two standard one-dimensional facts.  First, the
Fefferman--Stein decomposition \cite{FeffermanStein1972} says that
every mean-zero \(g\in\BMO(\T)\) has a representation
\begin{equation}
g=a+\mathcal Hb,\qquad
\|a\|_\infty+\|b\|_\infty
 \le C[g]_{\BMO(\T)}.
\label{eq:FS}
\end{equation}
Second, Hunt's endpoint estimate for ordinary and conjugate Fourier
sums \cite{Hunt1968,Oskolkov1974} gives absolute \(A,c>0\) such
that every \(u\in L^\infty(\T)\) satisfies
\begin{equation}
m\!\left\{x:\sup_{n\ge0}
  \max\{|S_nu(x)|,|\widetilde S_nu(x)|\}>t\right\}
\le A\exp(-c t/\|u\|_\infty),
\quad t>0.
\label{eq:hunt-endpoint}
\end{equation}
For \(u=0\), the left-hand side is zero.

The estimate for \(\mathcal Cg=\sup_n|S_ng|\) in
\cite[Remark 6.1]{Karagulyan2026} is an
\(L^\infty\!\to\BMO\) statement.

\begin{lemma}
\label{lem:exp}
There are absolute \(A,c>0\) such that, for every mean-zero
\(g\in\BMO(\T)\) and \(u>0\),
\begin{equation}
m\left\{x:\sup_{n\ge0}|g(x)-S_ng(x)|
     >u[g]_{\BMO(\T)}\right\}\le Ae^{-cu}.
\label{eq:remainder-tail}
\end{equation}
If \([g]_{\BMO}=0\), the maximal remainder vanishes almost
everywhere.
\end{lemma}

\begin{proof}
Let \(\sigma=[g]_{\BMO(\T)}\).  If \(\sigma=0\), the mean-zero
function \(g\) vanishes almost everywhere.  Otherwise choose
\(a,b\) in \eqref{eq:FS}.  Since \(S_n\) commutes with
\(\mathcal H\),
\[
\sup_n|S_ng|
\le\sup_n|S_na|+\sup_n|\widetilde S_nb|.
\]
Apply \eqref{eq:hunt-endpoint} to \(a\) and \(b\), and use
\(\|a\|_\infty+\|b\|_\infty\le C\sigma\).
This yields
\[
m\{x:\sup_n|S_ng(x)|>u\sigma\}\le A_1e^{-c_1u}.
\]
John--Nirenberg \cite{JohnNirenberg1961}, applied to the whole
circle and using \(g_\T=0\), likewise yields
\(m\{|g|>u\sigma\}\le A_2e^{-c_2u}\).
Finally \(\sup_n|g-S_ng|\le|g|+\sup_n|S_ng|\);
a union bound proves \eqref{eq:remainder-tail}.
\end{proof}

Let \(\mathcal M^{(1)}\) and \(\mathcal M^{(2)}\) be the
Hardy--Littlewood maximal operators along the indicated coordinates,
applied to the absolute value of their arguments.

\begin{lemma}\label{lem:dirichlet}
For \(h\in L^1(\T^2)\), \(m\ge0\), and \(\nu=1,2\),
\begin{equation}
|S_m^{(\nu)}h|
 \le C\log(m+2)\,\mathcal M^{(\nu)}h
\quad\text{almost everywhere.}
\label{eq:dirichlet}
\end{equation}
\end{lemma}

\begin{proof}
In the active coordinate the Dirichlet kernel obeys
\(|D_m(t)|\le C\min\{m+1,|t|^{-1}\}\).  Split the
convolution into \(|t|\le(m+1)^{-1}\) and dyadic annuli.
Each piece is bounded by a fixed multiple of the one-dimensional
maximal function, and there are \(O(\log(m+2))\) pieces.
\end{proof}

We also use, for \(1<s<\infty\),
\begin{equation}
\|\mathcal M^{(\nu)}H\|_{L^s(\T^2)}
 \le C_s\|H\|_{L^s(\T^2)}.
\label{eq:maximal-Ls}
\end{equation}

\section{Proof of Theorem~\ref{thm:main}}

Put \(\delta_q=(n_q+1)^{-1}\), and for \(\nu=1,2\) set
\(\varepsilon_q^{(\nu)}=C_J\Omega_\nu(f,\delta_q)\).
Choose \(P_q^{(\nu)}=P_{n_q}^{(\nu)}\) from
Lemma~\ref{lem:jackson} and write
\(g_q^{(\nu)}=f-P_q^{(\nu)}\).  Each active-coordinate
section of \(g_q^{(\nu)}\) has zero mean and BMO seminorm at most
\(\varepsilon_q^{(\nu)}\), outside a null set of passive coordinates.
Define
\[
\mathfrak B_q^{(\nu)}
 =\sup_{n_q\le n<n_{q+1}}|f-S_n^{(\nu)}f|.
\]
Since \(S_n^{(\nu)}P_q^{(\nu)}=P_q^{(\nu)}\) for \(n\ge n_q\),
we have
\[
\mathfrak B_q^{(\nu)}
 =\sup_{n_q\le n<n_{q+1}}
   |g_q^{(\nu)}-S_n^{(\nu)}g_q^{(\nu)}|.
\]
  Write
\(\varepsilon=\varepsilon_q^{(1)}\) and, for a passive coordinate
\(y\), put \(g_y=g_q^{(1)}(\cdot,y)\) and
\(\sigma_y=[g_y]_{\BMO(\T)}\).  For almost every \(y\),
\(g_y\) has mean zero and \(\sigma_y\le\varepsilon\).  Also,
\[
\mathfrak B_q^{(1)}(x,y)
 \le\sup_{n\ge0}|g_y(x)-S_ng_y(x)|.
\]
Suppose \(\varepsilon>0\).  If \(\sigma_y>0\), then
\(u\varepsilon\ge u\sigma_y\), so Lemma~\ref{lem:exp} yields
\[
\begin{aligned}
m\{x:\mathfrak B_q^{(1)}(x,y)>u\varepsilon\}
 &\le m\left\{x:\sup_{n\ge0}|g_y(x)-S_ng_y(x)|
            >u\sigma_y\right\}\\
 &\le Ae^{-cu},\qquad u>0.
\end{aligned}
\]
If \(\sigma_y=0\), the mean-zero section \(g_y\) vanishes
almost everywhere and the set on the left is empty.  Thus the
same bound holds for almost every \(y\).  Since \(m(\T)=1\),
Fubini gives
\[
\begin{aligned}
m_2\{\mathfrak B_q^{(1)}>u\varepsilon\}
 &=\int_\T m\{x:\mathfrak B_q^{(1)}(x,y)>u\varepsilon\}\,dy\\
 &\le Ae^{-cu},\qquad u>0.
\end{aligned}
\]
For \(\nu=2\), fix \(x\) instead and integrate over \(x\).
We have therefore proved, for both \(\nu=1,2\) whenever
\(\varepsilon_q^{(\nu)}>0\), that
\begin{equation}
m_2\{\mathfrak B_q^{(\nu)}>u\varepsilon_q^{(\nu)}\}
 \le Ae^{-cu},\qquad u>0.
\label{eq:block-tail}
\end{equation}
Finally, if
\(\varepsilon_q^{(\nu)}=0\), almost every active-coordinate
section of \(g_q^{(\nu)}\) has BMO seminorm zero and mean zero.
Such a section is zero almost everywhere, so
\(\mathfrak B_q^{(\nu)}=0\) almost everywhere by Fubini.

We construct \(\lambda_q\ge1\) with
\begin{equation}
\sum_{q\ge1}e^{-\lambda_q}<\infty,
\qquad a_q\lambda_q\longrightarrow0.
\label{eq:lambda-goal}
\end{equation}
Indeed, \eqref{eq:new-condition} implies \(a_q\to0\).
For \(k\ge1\), choose \(1\le Q_1<Q_2<\cdots\)
so that, for \(q\ge Q_k\),
\(a_q\le1/k\) and
\[
\sum_{q\ge Q_k}\exp[-1/(ka_q)]\le2^{-k}.
\]

For \(Q_k\le q<Q_{k+1}\), set
\(\lambda_q=(ka_q)^{-1}\) when \(a_q>0\), and
\(\lambda_q=q\) otherwise.  For the finitely many \(q<Q_1\),
put \(\lambda_q=q\).  The sums over successive blocks are at
most \(2^{-k}\), while the zero cases contribute at most
\(\sum_qe^{-q}\).  This proves the first part of
\eqref{eq:lambda-goal}; on a nonzero block
\(a_q\lambda_q=1/k\), and in zero cases the product is
zero.  Thus the second part holds as well.
In particular,
\begin{equation}
2^q\varepsilon_q^{(\nu)}\lambda_q
 \le C_Ja_q\lambda_q\longrightarrow0.
\label{eq:small-blocks}
\end{equation}

Choose \(K\) with \(cK\ge2\), where \(c\) is from
\eqref{eq:block-tail}, and set
\[
C_\nu=\max\left\{1,\sup_{q\ge1}
\frac{\mathfrak B_q^{(\nu)}}
{K\varepsilon_q^{(\nu)}\lambda_q}\right\}.
\]
Declare a quotient with \(\varepsilon_q^{(\nu)}=0\) to be
zero, on a common full-measure set where the corresponding block
vanishes.  For \(t>1\), \eqref{eq:block-tail} gives
\begin{align*}
m_2\{C_\nu>t\}
 &\le A\sum_{q\ge1}e^{-cKt\lambda_q}\\
 &\le Ae^{-t}\sum_{q\ge1}e^{-\lambda_q},
\end{align*}
because \(\lambda_q\ge1\) and \(cK\ge2\).
Consequently \(C_\nu\in L^s(\T^2)\) for every finite
\(s\ge1\), and
\begin{equation}
\mathfrak B_q^{(\nu)}
 \le K\varepsilon_q^{(\nu)}\lambda_q C_\nu
\quad\text{almost everywhere for every }q.
\label{eq:block-pointwise}
\end{equation}
By \eqref{eq:maximal-Ls},
\(\mathcal M^{(1)}C_2\) and \(\mathcal M^{(2)}C_1\)
are finite almost everywhere.

Suppose \(M\le N\), and choose \(q\) with
\(n_q\le N<n_{q+1}\).  Then
\[
f-S_{M,N}f
 =(f-S_M^{(1)}f)
  +S_M^{(1)}(f-S_N^{(2)}f).
\]
For almost every \(y\), \(f(\cdot,y)\in\BMO(\T)\)
and has a finite mean.  John--Nirenberg gives
\(f(\cdot,y)\in L^2(\T)\).  Carleson's theorem \cite{Carleson1966} and
Fubini therefore imply
\begin{equation}
S_M^{(1)}f(x,y)\longrightarrow f(x,y)
\quad\text{for almost every }(x,y).
\label{eq:first-coordinate}
\end{equation}
Using Lemma~\ref{lem:dirichlet},
\eqref{eq:block-pointwise}, and \(M\le N<n_{q+1}\),
we obtain, on the same full-measure set,
\begin{align*}
|S_M^{(1)}(f-S_N^{(2)}f)|
 &\le C\log(M+2)
       \mathcal M^{(1)}(f-S_N^{(2)}f)\\
 &\le C'2^q\varepsilon_q^{(2)}\lambda_q
       \mathcal M^{(1)}C_2
 \longrightarrow0
\end{align*}
by \eqref{eq:small-blocks}, as \(N\to\infty\).
When \(N<M\), reverse the roles of the two coordinates.
Together with \eqref{eq:first-coordinate}, this proves the
Pringsheim limit.

\section{Proof of Theorem~\ref{thm:example}}

For \(k\ge1\), set
\[
R_k=2^{2^{k+2}},\qquad
N_k=2^{2^{R_k}},\qquad
A_k=\frac1{2^{R_k}\log R_k}.
\]
Here \(R_{k+1}=R_k^2\), \(N_k=n_{R_k}\), and
\begin{equation}
\frac{A_kN_k}{A_{k+1}N_{k+1}}<\frac12,\quad
\frac{A_{k+1}}{A_k}<\frac12,\quad
A_{k-1}<N_{k-1}A_k\quad(k\ge2).
\label{eq:scale-inequalities}
\end{equation}

Define
\[
g(x)=\sum_{k=1}^\infty A_k\sin(2\pi N_kx),
\qquad F(x,y)=g(x)+g(y).
\]
The series converges uniformly because \(\sum_kA_k<\infty\),
so \(F\in C(\T^2)\).

Put \(d_q=2^{-2^q}=1/n_q\).  If \(k\ge2\) and
\(R_{k-1}<q\le R_k\), then \(d_q\le N_{k-1}^{-2}\).
For every \(|h|\le d_q\), \eqref{eq:scale-inequalities}
gives
\begin{align*}
|g(x+h)-g(x)|
 &\le2\pi d_q\sum_{\ell<k}A_\ell N_\ell
       +2\sum_{\ell\ge k}A_\ell\\
 &\le4\pi\,\frac{A_{k-1}}{N_{k-1}}+4A_k
 \le C A_k.
\end{align*}
 Hence
\begin{equation}
\omega(F,d_q)_C\le C' A_k
\qquad(R_{k-1}<q\le R_k).
\label{eq:example-upper}
\end{equation}
Since \(\omega(F,(n_q+1)^{-1})_C\le\omega(F,d_q)_C\),
the denominator \(b_q=2^q\omega(F,(n_q+1)^{-1})_C\)
in \eqref{eq:example-condition} satisfies
\[
b_q\le C'2^q A_k
 =\frac{C'2^{q-R_k}}{\log R_k}.
\]
Thus, for every \(\eta>0\), the sum over one block is
\begin{align*}
\sum_{R_{k-1}<q\le R_k}e^{-\eta/b_q}
 &\le\sum_{j=0}^\infty
 \exp\!\left[-\frac{\eta}{C'}\,2^j\log R_k\right]\\
 &\le C_\eta R_k^{-\eta/C'}.
\end{align*}
The last expression is summable in \(k\), since \(R_k\)
grows doubly exponentially.  The finitely many initial
terms cause no difficulty.  This proves
\eqref{eq:example-condition}.

Now set \(h_k=(4N_k)^{-1}<d_{R_k}\).
For \(\ell>k\), \(N_\ell/N_k\) is divisible by \(4\), so
\(\sin(2\pi N_\ell h_k)=0\).  For \(\ell<k\),
\(0<2\pi N_\ell h_k<\pi/2\), so those sine terms are
nonnegative.  It follows that
\[
\omega(F,d_{R_k})_C
\ge |F(h_k,0)-F(0,0)|\ge A_k.
\]
With \(r_k=2^{R_k}\), we have
\(d_{R_k}=2^{-r_k}\) and
\[
r_k\log\log r_k\;\omega(F,2^{-r_k})_C
\ge\frac{\log(R_k\log2)}{\log R_k}
\longrightarrow1.
\]
This proves \eqref{eq:example-failure}.  In fact, the
specific \(F\) has uniformly convergent rectangular Fourier
sums, since its Fourier coefficients are absolutely summable.

\end{document}